
\input amssym.def
\input amssym

\documentstyle{amsppt}
\magnification=1200

\hoffset 0,7truecm
\voffset 1truecm
\hsize 15truecm
\vsize 22truecm 

\catcode`\@=11
\long\def\thanks#1\endthanks{%
  \ifx\thethanks@\empty \gdef\thethanks@{%
     \tenpoint#1} 
     \else \expandafter\gdef\expandafter\thethanks@\expandafter{%
     \thethanks@\endgraf#1}%
  \fi}
\def\keywords{\let\savedef@\keywords
  \def\keywords##1\endkeywords{\let\keywords\savedef@
  \toks@{\def\usualspace{{\it\enspace}}\tenpoint} 
  \toks@@{##1\unskip.}%
  \edef\thekeywords@{\the\toks@\frills@{%
  {\noexpand\it 
  Key words and phrases.\noexpand\enspace}}\the\toks@@}}%
  \nofrillscheck\keywords}
\def\makefootnote@#1#2{\insert\footins
 {\interlinepenalty\interfootnotelinepenalty
 \tenpoint\splittopskip\ht\strutbox\splitmaxdepth\dp\strutbox 
 \floatingpenalty\@MM\leftskip\z@\rightskip\z@
 \spaceskip\z@\xspaceskip\z@
 \leavevmode{#1}\footstrut\ignorespaces#2\unskip\lower\dp\strutbox
 \vbox to\dp\strutbox{}}}
\def\subjclass{\let\savedef@\subjclass
  \def\subjclass##1\endsubjclass{\let\subjclass\savedef@
  \toks@{\def\usualspace{{\rm\enspace}}\tenpoint}
  \toks@@{##1\unskip.}%
  \edef\thesubjclass@{\the\toks@\frills@{{%
  \noexpand\rm1991 {\noexpand\it Mathematics
  Subject Classification}.\noexpand\enspace}}%
  \the\toks@@}}%
\nofrillscheck\subjclass}
\def\abstract{\let\savedef@\abstract
 \def\abstract{\let\abstract\savedef@
  \setbox\abstractbox@\vbox\bgroup\noindent$$\vbox\bgroup
  \def\envir@end{\endabstract}\advance\hsize-2\indenti
  \def\usualspace{\enspace}\tenpoint \noindent 
  \frills@{{\smc Abstract.\enspace}}}%
 \nofrillscheck\abstract}

\outer\def\endtopmatter{\add@missing\endabstract
 \edef\next{\the\leftheadtoks}\ifx\next\empty
  \expandafter\leftheadtext\expandafter{\the\rightheadtoks}\fi
 \ifmonograph@\else
   \ifx\thesubjclass@\empty\else \makefootnote@{}%
        {\thesubjclass@}\fi
   \ifx\thekeywords@\empty\else \makefootnote@{}%
        {\thekeywords@}\fi  
   \ifx\thethanks@\empty\else \makefootnote@{}%
        {\thethanks@}\fi   
 \fi
  \pretitle
  \begingroup 
  \ifmonograph@ \topskip7pc \else \topskip4pc \fi
  \box\titlebox@
  \endgroup
  \preauthor
  \ifvoid\authorbox@\else \vskip2.5pcplus1pc\unvbox\authorbox@\fi
  \preaffil
  \ifvoid\affilbox@\else \vskip1pcplus.5pc\unvbox\affilbox@\fi
  \predate
  \ifx\thedate@\empty\else
       \vskip1pcplus.5pc\line{\hfil\thedate@\hfil}\fi
  \preabstract
  \ifvoid\abstractbox@\else
       \vskip1.5pcplus.5pc\unvbox\abstractbox@ \fi
  \ifvoid\tocbox@\else\vskip1.5pcplus.5pc\unvbox\tocbox@\fi
  \prepaper
  \vskip2pcplus1pc\relax}
\def\foliofont@{\tenrm} 
\def\headlinefont@{\tenpoint} 


\def\refstyle#1{\uppercase{%
  \if#1A\relax \def\keyformat##1{[##1]\enspace\hfil}%
  \else\if#1B\relax 
      \def\keyformat##1{\aftergroup\kern
           aftergroup-\aftergroup\refindentwd}%
      \refindentwd\parindent 
  \else\if#1C\relax 
      \def\keyformat##1{\hfil##1.\enspace}%
  \else\if#1D\relax 
      \def\keyformat##1{\hfil\llap{{##1}\enspace}} 
  \fi\fi\fi\fi}
}
\refstyle{D}   
\def\refsfont@{\tenpoint} 
\def\email{\let\savedef@\email
  \def\email##1\endemail{\let\email\savedef@
  \toks@{\def\usualspace{{\it\enspace}}\endgraf\indent\tenpoint}%
  \toks@@{\tt ##1\par}%
  \expandafter\xdef\csname email\number\addresscount@\endcsname
  {\the\toks@\frills@{{\noexpand\smc E-mail address\noexpand\/}:%
     \noexpand\enspace}\the\toks@@}}%
  \nofrillscheck\email}

\def\proof{\noindent{\bf Proof.\ }}


\newcount\secno

\secno=1

\topmatter

\title
Galois theory for a class of complete\\
modular lattices
\endtitle

\author
Alexandre A.$\,$Panin
\endauthor

\affil
{\it Department of Mathematics and Mechanics\\ 
St.Petersburg State University\\ 
2 Bibliotechnaya square,\\
St.Petersburg 198904, Russia}
\endaffil


\address
\endaddress
\email
alex@ap2707.spb.edu
\endemail


\thanks
\endthanks

\abstract
We construct Galois theory for sublattices of certain
complete modular lattices and their automorphism groups.
A well--known description of the intermediate subgroups
of the general linear group over a semilocal ring containing
the group of diagonal matrices, due to Z.I.Borewicz and
N.A.Vavilov, can be obtained as a consequence of this theory.
Bibliography: 3 titles.
\endabstract

\endtopmatter

\document

\heading
Introduction
\endheading

We generalize here the results of [PY], [S]. Namely, Galois theory for a
class of complete modular lattices is constructed.

By an automorphism of a complete lattice we mean hereafter a bijective
mapping of the lattice onto itself which commutes with the supremum and
infimum of every subset of the lattice.
Other notions and definitions are introduced in [PY].

\heading
Formulation of the main results
\endheading

Let $L$ be a complete modular lattice, $L_0$ its finite sublattice, which is
a Boolean algebra, $G$ a subgroup of the group of all automorphisms of
the lattice $L$, $H=G(L_0)$.

Let $e_1,e_2,\ldots,e_n$ be the atoms of $L_0$. We consider a number of
additional conditions (it is supposed that, unless otherwise stated,
the indices are changing from $1$ to $n$):

\smallbreak

$1^0.\ \ 0_{L_0}=0_L,\ 1_{L_0}=1_L$.

\smallskip

$2^0$. If $f\in G$ and $f(e_i)+\widehat{e}_j=1$ for some $i,j$, then
$f(e_i)\cdot\widehat{e}_j=0$.

\smallskip

$3^0$. If $a\in G$ and $[a(e_i)]_i=e_i$ for some $i$, then there exists
$h\in H_i$ such that $[ha(x_i)]_i=[ah(x_i)]_i=x_i$ for every
$x_i\leqslant e_i$.

\smallskip

$4^0$. There exists $h\in H_t\cap G(\overline{L}_0)$ such that
$[aha^{-1}(x_i)]_r = [a([a^{-1}(x_i)]_t)]_r$ for every
$a\in G,\ r\neq i,\ x_i\leqslant e_i$.

\smallskip

$5^0$. Let $u\in\overline{L}_0,\ u\geqslant e_i$ for some
$i$; let $g\in G,\ [g(u)]_i=e_i$. Then there exists $t\in G$ such that:

$1)\ [gt(e_i)]_i=e_i$,

$2)\ t(e_s)=e_s$ for $s\neq i$,

$3)\ [t(e_i)]_j\leqslant [u]_j$.

\smallskip

$6^0$. If $f,g\in G$ and $[f(e_i)]_j\leqslant [g(e_i)]_j$ for
some $i,j$, then $[f(x)]_j\leqslant [g(x)]_j$ for every
$x\in L_0^{'},\ x\leqslant e_i$.

\smallskip

$7^0$. If $u\leqslant e_j$ for some $j$, then for every $i\neq j$
there exist $y_{\alpha}\leqslant e_j, \alpha\in I$ such that
$u=\sum\limits_{\alpha\in I}y_{\alpha}$ and
$H_{ij}(y_{\alpha})\neq\varnothing$.

\smallskip

$8^0$. If $x=[f(e_i)]_j$ for some $f\in G,\ i\neq j$, then there exists
$g\in H_{ij}(x)$ such that $[g(u)]_j=[f(u)]_j$ for every $u\leqslant\ e_i$.

\smallskip

$9^0$. If $w\in L,\ [w]=(0,\ldots,e_i,\ldots,x,\ldots,0)$, where
$w\cdot e_j=0$ and $H_{ij}(x)\neq\varnothing$,
then there exists $t\in H_{ij}(x)$ such that $t(w)=e_i$.

\smallskip

$10^0$. Let $a_{\alpha}\in H_{ij}(x_{\alpha})$ and
$y\leqslant\sum\limits_{\alpha\in I}x_{\alpha}$ be such that
$H_{ij}(y)\neq\varnothing$. Then
$H_{ij}(y)\subseteq \langle H,a_{\alpha}:\ \alpha\in I\rangle$.

\smallskip

$11^0$. If $a\in G$, then for every $i\neq j$ and every $h\in H_t$ the set
$H_{ij}([aha^{-1}(e_i)]_j)\cap \langle a,H\rangle$ is not empty.

\smallskip

$12^0.\ \ L_0^{'}\subseteq \overline{L}_0$.

\smallskip

$13^0$. For every $x\leqslant e_i,\ x\neq e_i$ there exists a coatom
$y\leqslant e_i$ such that $x\leqslant y$.

\smallskip

We denote $w_i=\prod\limits_{x\,is\,a\,coatom\,in\,e_i}x\quad$ for every $i$
and $w=\sum\limits_{i=1}^{n}w_i$.

\smallskip

$14^0$. The lattice $L^w=\{x\in L:\ x\geqslant w\}$ is of finite length.

\smallskip

We denote $G^w=\{g\in G:\ g$ is identical on $L^w\}$.

\smallskip

$15^0$. Let $t\in H_{ij}(x),\ x\leqslant w_j$. Then there exists
$h\in H$ such that $th\in G^w$.

\smallskip

$16^0$. The lattice $L^w\cap\overline{L}_0(H)$ is finite.

\proclaim
{Theorem 1}
Assuming that the conditions $1^0-12^0$ are fulfilled,
for every subgroup $F\geqslant H$ of group $G$:
\smallskip\indent
(i)$\ \sigma=\sigma(F)$ is a net collection in $L_0^{'}$;
\smallskip\indent
(ii)$\ G(K_{\sigma})\trianglelefteq F$;
\smallskip\indent
(iii) if $M$ is a sublattice of $L_0^{'}$ such that
$G(M)\trianglelefteq F$, then $G(M)=G(K_{\sigma})$.
\endproclaim

\proclaim 
{Theorem 2}
Let $\tau=(\tau_{ij})$ be a net collection in
$L_0^{'},\ g\in G$. Provided that the conditions $1^0-11^0,\ 13^0-16^0$
are fulfilled, we have:
\smallskip\indent
(i) if $[g(e_i)]_j\leqslant \tau_{ij}$ for every $i,j$, then
$g\in G(K_{\tau})$;
\smallskip\indent
(ii) the index of $G(K_{\tau})$ in its normalizer is finite.
\endproclaim

\heading
Proof of the main results
\endheading

Proof of the parts (i)--(iii) of Theorem 1 is analogous to the proof of the
corresponding assertions of [PY]. We note that instead of
properties of the dimension function on $L$ used in [PY], one must apply
the modularity law and the following trivial statement:

\smallskip

If $A\subseteq G$, where $A^{-1}=A$ and $a(x)\leqslant x$ for every $a\in A$
and some $x\in L$, then $a(x)=x$ for every $a\in A$.

\smallskip

Proof of Theorem 2 will be presented below.

\proclaim 
{Lemma 1}
For every $f\in G\ f(w)=w$.
\endproclaim

\proof It is sufficient to check $[f(w_i)]_j\leqslant w_j$ for every $i,j$.

It is clear that $w\in L_0^{'}$, therefore for $i=j$ it is just the
condition $6^0$.

Let $i\neq j,\ g\in H_{ij}(e_j)$. By the condition $9^0$ there exists
$t\in H_{ji}(e_i)$ such that $tg(e_i)=e_j$. Consider an arbitrary
coatom $y$ in $e_j$. Then $(tg)^{-1}(y)$ is a coatom in $e_i$, therefore
$tg(w_i)\leqslant y$, whence $[g(w_i)]_j\leqslant w_j$. Now it remains
to apply the condition~$6^0$.

\smallskip
Suppose $\tau=(\tau_{ij})$ is a net collection in $L_0^{'}$.

\proclaim
{Lemma 2}
Let $\rho_{ij}=\tau_{ij}+w_j$. Then:
\smallskip\indent
(i)$\ \ \rho=(\rho_{ij})$ is a net collection in $L_0^{'}$;
\smallskip\indent
(ii)$\ \ G(K_{\rho})=G(K_{\tau})\cdot G^w$.
\endproclaim

\proof Let $i,j,k$ be pairwise distinct, and let
$g\in H_{ij}(x),\ x\leqslant\tau_{ij}+w_j$.

By the conditions $7^0$ and
$10^0\ g\in \langle H,H_{ij}(y),H_{ij}(z):\ y
\leqslant\tau_{ij},\ z\leqslant w_j\rangle$.

If $f\in H_{ij}(y)$, where $y\leqslant\tau_{ij}$, then
$[f(\tau_{ki}+w_i)]_j\leqslant\tau_{kj}+w_j$ by Lemma 1. If $f\in H_{ij}(z)$,
where $z\leqslant w_j$, then $[f(\tau_{ki}+w_i)]_j\leqslant w_j$.

(i) Apply the condition $8^0$.

(ii) The inclusion $\supseteq$ is trivial. Further, by the condition $15^0$ and
by Theorem~7.2~[PY] $G(K_{\rho})\subseteq \langle G(K_{\tau}),G^w\rangle$. It
remains to note that $G^w\trianglelefteq G$.

\smallskip

\noindent {\bf Proof of Theorem 2(i).}

A. First, suppose that $\tau_{ij}\leqslant w_j$ for every $i\neq j$. Since
$1=\sum\limits_{i=1}^{n}g(e_i)$, we have $e_1=[g(e_1)]_1+w_1$. By the condition
$13^0\ \ [g(e_1)]_1=e_1$. Repeating the proof of Theorem~7.2~[PY], we obtain
$g\in \langle H,\ H_{ij}(x):\ x\leqslant \tau_{ij}\rangle\leqslant
G(K_{\tau})$.

B. General case. We put $x_i=\sum\limits_{j=1}^{n}\tau_{ij}$.
By the definition of a net collection we have $g(x_i)\leqslant x_i$.
Further, it follows from Lemma 1 that $g(x_i+w)\leqslant x_i+w$, and since
the restriction of $g$ to $L^w$ is an automorphism of this lattice, we have
$g(x_i+w)=x_i+w$ by the condition $14^0$. Thus $g\in G(K_{\rho})$.

It follows from Lemma 2 that $g=g_1g_2$, where
$g_1\in G(K_{\tau}),\ g_2\in G^w$. Since
$g_2(x_i)\leqslant x_i$, then $[g_2(e_i)]_j\leqslant\tau_{ij}\cdot w_j$
for every $i\neq j$. We have already proved in the part A that
$[g_2^{-1}(e_i)]_j\leqslant\tau_{ij}\cdot w_j$
for every $i\neq j$, therefore $g_2\in G(K_{\tau})$.

\smallskip

Proof of the part (ii) is conceptually identical with the
constructions of \S$\,$7 of the article [BV].

\smallskip

As in [PY], a complete description of subgroups of the general linear group
over a semilocal ring (whose fields of residues have at least
seven elements, see [BV]), containing the group of diagonal matrices,
can be deduced from Theorems 1\break and 2.

\heading
Acknowledgements
\endheading

The author is grateful to Professor Anatoly V.$\,$Yakovlev for the
formulation of the problem and useful discussions.


\frenchspacing
\Refs

\ref \key\bf{{[BV]}} \by Borewicz Z.I., Vavilov N.A.
\paper Subgroups of the full linear group over a semilocal
ring containing the group of diagonal matrices
\jour Proc. Steklov Inst. Math.
\yr 1980 \issue 4 \pages 41--54 \endref

\ref \key\bf{{[PY]}} \by Panin A.A., Yakovlev A.V.
\paper Galois theory for a class of modular lattices
\jour Zap. Nauchn. Semin. POMI
\yr 1997 \vol 236 \pages 133--148 (In Russian, English transl.
to appear in J.~Math. Sci.)
\endref

\ref \key{\bf{[S]}} \by Simonian A.Z.
\paper Galois theory for modular lattices
\inbook Ph.D. thesis, St.Pe\-ter\-sburg State University 
\yr 1992 \pages 1--73 (In Russian)
\endref 

\endRefs
\enddocument